\documentclass{amsart}

\usepackage{amsfonts}
\usepackage{amsmath}
\usepackage[all]{xy}

\def\xypic{\hbox{\rm\Xy-pic}}

\DeclareFontFamily{OMS}{rsfs}{\skewchar\font'60}
\DeclareFontShape{OMS}{rsfs}{m}{n}{<-5>rsfs5 <5-7>rsfs7 <7->rsfs10 }{}
\DeclareSymbolFont{rsfs}{OMS}{rsfs}{m}{n}
\DeclareSymbolFontAlphabet{\scr}{rsfs}

\DeclareMathOperator{\Id}{Id}
\DeclareMathOperator{\id}{id}
\DeclareMathOperator{\Ext}{Ext}

\begin{document}


\author {Anna Marie Bohmann and J. P. May}



\address{Department of Mathematics, Northwestern University\\ 2033 Sheridan Road, Evanston IL 60208\\[5pt] Department of Mathematics, University of Chicago\\5734 S. University Avenue, Chicago IL 60637\\}

\title{A presheaf interpretation of the generalized Freyd conjecture}







\def\xypic{\hbox{\rm\Xy-pic}}

\newcommand{\sA}{\scr{A}}
\newcommand{\sB}{\scr{B}}
\newcommand{\sC}{\scr{C}}
\newcommand{\sD}{\scr{D}}
\newcommand{\sM}{\scr{M}}
\newcommand{\sP}{\scr{P}}
\newcommand{\sS}{\scr{S}}
\newcommand{\sT}{\scr{T}}

\newcommand{\bF}{\mathbb{F}}
\newcommand{\bP}{\mathbb{P}}
\newcommand{\bU}{\mathbb{U}}
\newcommand{\bY}{\mathbb{Y}}
\newcommand{\bZ}{\mathbb{Z}}

\newcommand{\epz}{\varepsilon}
\newcommand{\io}{\iota}
\newcommand{\SI}{\Sigma}

\newcommand{\com}{\circ}     
\newcommand{\iso}{\cong}     
\newcommand{\rtarr}{\rightarrow}
\newtheorem{thm}{Theorem}

\newtheorem{prop}{Proposition}
\newtheorem{lem}{Lemma}
\newtheorem{cor}{Corollary}
\newtheorem{conj}{Conjecture}

\theoremstyle{definition}
\newtheorem{defn}[thm]{Definition}
\newtheorem{exmp}[thm]{Example}

\theoremstyle{remark}
\newtheorem{rem}[thm]{Remark}

\let\pf\proof
\let\epf\endproof


\maketitle

\begin{abstract}
We give a generalized version of the Freyd conjecture and a way
to think about a possible proof.  The essential point is to describe
an elementary formal reduction of the question that holds in any triangulated category. 
There are no new results, but
at least one known example drops out very easily.
\end{abstract}


In algebraic topology, the generating hypothesis, or Freyd conjecture, is a long-standing conjecture about the structure of the stable homotopy category.  It was initially formulated in 1965 \cite{F} and remains open to this day.  Because the original conjecture has proved difficult to analyze, recent work has turned to studying similar conjectures in categories that share many properties with the stable homotopy category in hopes of further understanding the types of categories in which such a conclusion holds.  In this note, we state a version of the generating hypothesis for an arbitrary triangulated category and analyze conditions under which this hypothesis holds.  We emphasize that we impose no additional conditions on our triangulated categories, so that our results show which 
{\em formal} properties of a category imply the Freyd conjecture.  Therefore our results give conceptual insight into the kind of category in which the generating hypothesis can be expected to hold.

\section{The generalized Freyd conjecture}

Let $\sT$ be a triangulated category and
write $[X,Y]$ for the abelian group of maps $X\rtarr Y$ in $\sT$. Let 
$\sB$ be a (small) full subcategory of $\sT$ closed under its translation
(or shift) functor $\SI$ and let $\sC$ be the thick full subcategory of $\sT$
that is generated by $\sB$; write $\io\colon \sB\rtarr \sC$ for the inclusion.
For emphasis, we often write $\sB(X,Y) = [X,Y]$ when $X,Y\in\sB$ and
$\sC(X,Y) = [X,Y]$ when $X,Y\in\sC$.  The category $\sB$ is pre-additive 
(enriched over $\sA\!b$), and $\sC$ is additive (has biproducts).  Let 
$\sP\sB$ and $\sP\sC$ denote the categories of abelian presheaves defined 
on $\sB$ and $\sC$. They consist of the additive functors from $\sB^{op}$ 
or $\sC^{op}$ to $\sA\!b$ and the additive natural transformations between them.

\begin{defn}\label{defn:Ffunctor} Define the Freyd functor $\bF\colon \sT\rtarr \sP\sB$
by sending an object $X$ to the functor $\bF X$ specified on objects 
and morphisms of $\sB$ by $\bF X(-) = [-,X]$ and sending a map 
$f\colon X\rtarr Y$ to the natural transformation $f_* =[-,f]$.  
Define $\bY\colon \sT\rtarr \sP\sC$ similarly.  We are only interested in the 
restrictions of $\bF$ and $\bY$ to $\sC$, and then $\bY\colon \sC\rtarr \sP\sC$ 
is the standard Yoneda embedding.
\end{defn}

\begin{conj}[The generalized Freyd conjecture]\label{FC}  The functor
$\bF\colon \sC\rtarr \sP\sB$ is faithful.  Equivalently,  $\bF f =0$ if 
and only if $f=0$.  We then say that the Freyd conjecture holds for the 
pair $(\sC,\sB)$.
\end{conj}

We hasten to add that the conjecture is false in this generality.
Additional hypotheses are needed, but none will be relevant
to our formal analysis.  For example, we might as well assume 
that a map $f\colon X\rtarr Y$ in $\sC$ is an isomorphism if and only 
$f_*\colon \bF X\rtarr \bF Y$ is an isomorphism, as holds in the
stable homotopy category. This condition is necessary but not sufficient 
for the Freyd conjecture to hold. Indeed, if $Z$ is the third term in an 
exact triangle with one map $f$, then $f$ is an isomorphism if and only 
if $Z=0$ and $\bF f$ is an isomorphism if and only if $\bF Z=0$. If 
$\bF Z = 0$ and the Freyd conjecture is true, then $Z=0$ since its
identity map is the zero map.

\begin{exmp}\label{FC1} Take $\sT = \text{Ho}\sS$ to be the 
stable homotopy category.  Let $\sB$ consist of the sphere spectra
$S^n$ for integers $n$.  Then $\sC$ is the homotopy category of finite 
CW spectra.  Freyd \cite{F} conjectured that a map $f$ in $\sC$ is zero if and
only if it induces the zero homomorphism $f_*\colon \pi_*(X)\rtarr \pi_*(Y)$.
By the following observation, this is a special case of our
Conjecture \ref{FC}.
\end{exmp}

\begin{lem}\label{FO} In Example \ref{FC1}, the category $\sP\sB$ is isomorphic 
to the category $\sM$ of right modules over the ring $\pi_*$ of stable 
homotopy groups of spheres.  Under this isomorphism, the Freyd functor $\bF$ 
coincides with the stable homotopy group functor $\pi_*\colon \text{Ho}\sS\rtarr \sM$.
\end{lem}
\begin{proof}  We have $\sB(S^m,S^n) = \pi_m(S^n)\iso \pi_{m-n}$. For $T\in \sP\sB$, 
let $T_n = T(S^n)$.  The additive contravariant functor $T$ gives homomorphisms
\[ T\colon \pi_{m-n} \iso \sB(S^m,S^n)\rtarr \sA\!b(T_n,T_m).\]
By adjunction, these give homomorphisms $T_n \otimes \pi_{m-n} \rtarr T_m$. 
The functoriality gives the formula $(tx)y = t(xy)$.  No signs appear since
we are taking right modules, as is dictated categorically by contravariance.
Conversely, given a right $\pi_*$-module $M$, define $T(S^n)=M_n$ and, for 
$t\in M_n$ and $x\in \pi_{m-n}\iso \sB(S^m,S^n)$  define $T(x)(t) = tx\in M_m$. 
The module axioms ensure that $T$ is a functor.  This specifies the isomorphism
of categories, and the consistency of $\bF$ and $\pi_*$ is clear.
\end{proof} 

The following example originally prompted us to take a presheaf perspective 
on the Freyd conjecture.

\begin{exmp}\label{FCG} Let $G$ be a compact Lie group and take $\sT = \text{Ho}G\sS$ 
to be the equivariant stable homotopy category.  Let $G\sB$ consist of the orbit 
$G$-spectra 
$$S^n[G/H]\equiv \SI^n\SI^{\infty}(G/H)_+$$ 
for integers $n$ and closed 
subgroups $H$ of $G$.  The thick subcategory $G\sC$ generated by $G\sB$ is the category 
of retracts of finite $G$-CW spectra.  The equivariant version of the Freyd conjecture 
asserts that a map $f$ in $G\sC$ is zero if and only if it induces the zero homomorphism 
$f_*\colon \pi_*^H(X)\rtarr \pi_*^H(Y)$ for all $H$, where 
\[ \pi_n^H(X) \equiv \pi_n(X^H) \iso [S^n[G/H],X]_G. \]
Again, this is a special case of our Conjecture \ref{FC}.  This example is the focus of \cite{B}.
\end{exmp}

\begin{rem}\label{MacF}  The full subcategory $\sB_0$ of $\sB$ whose objects are the 
$S[G/H]\equiv S^0[G/H]$ is called the Burnside category.  A Mackey functor (or $G$-Mackey 
functor) $M$ is by definition an object of $\sP G\sB_0$.  When $G$ is finite, this agrees 
with the more usual algebraic definition (\cite[V\S9]{LMS} or \cite[IX\S4, XIX\S3]{EHCT}).
A map of Mackey functors is a natural transformation, that is, a map of presheaves. 
\end{rem}

\begin{defn} The graded Burnside category $G\pi_*$ has objects 
the $S[G/H]$.  Its abelian group of maps of degree $n$ from $S[G/H]$ to
$S[G/J]$ is $$\pi_n^H(S[G/J]) = [S^n[G/H],S[G/J]]_G.$$  Composition is induced by 
suspension and composition in $G\sB$ in the evident fashion. Define a right 
$G\pi_*$-module $M$ to be a graded presheaf, that is, a contravariant functor 
$G\pi_*\rtarr \sA\! b_*$, where $\sA\!b_*$ is the category of graded abelian groups.
\end{defn}

\begin{lem}\label{FOG} In Example \ref{FCG}, the category $\sP G\sB$ is isomorphic 
to the category $G\sM$ of right $G\pi_*$-modules.  Under this isomorphism, the 
Freyd functor $\bF$ coincides with the equivariant stable homotopy group functor 
$\pi_*^{(-)}\colon \text{Ho}G\sS\rtarr G\sM$.
\end{lem}

The proof is the same as that of Example \ref{FC1}. In this case, the presheaf formulation 
of the Freyd conjecture appears more natural to us than the equivalent homotopy group reformulation. 

\begin{exmp} In the rational equivariant stable homotopy category, the Freyd 
conjecture is true if $G$ is finite \cite{GM} and is false if $G = S^1$ \cite{B}.
\end{exmp}

Here is another example where the presheaf formulation, and the use of many objects and
their shifts rather than a single object and its shifts in $\sB$, may be more natural than the formulation of the Freyd 
conjecture studied so far.

\begin{exmp} In several recent papers, Benson, Carlson, Chebolu, Christensen and Min\'a\v c \cite{BCCM, CCM2008, CCM2007} study the Freyd conjecture in the stable module category of a finite group $G$ over a field $k$ whose characteristic divides the order of $G$.  This is a triangulated category $\mathrm{StMod}_{kG}$ obtained from the category of $kG$--modules by modding out by maps that factor through a projective module.  They restrict to the thick subcategory generated by the trivial representation $k$ and ask whether the Tate cohomology functor is faithful on this subcategory.  Since the Tate cohomology of a module $M$ is given by maps out of $k$ in the stable module category,
\[ \hat{H}^i(G,M)=\mathrm{Stmod}_{kG}(\Omega^i k,M),\] their formulation is 
equivalent to our Conjecture \ref{FC} with $\sB=\{\Omega^i k\}$.  In this context, they prove that the Freyd conjecture holds in the stable module category of $kG$--modules if and only if the 
$p$--Sylow subgroup of $G$ is either $C_2$ or $C_3$, where $p$ is the characteristic of $k$ \cite{CCM2008}.  Their proof in fact shows that this variant of the Freyd conjecture holds if and only if the thick subcategory generated by $k$ consists of finite direct sums of suspensions of $k$.

For $p$--groups, the trivial module is the only irreducible module over $kG$, and the thick subcategory generated by $k$ is the subcategory of compact objects in the stable module category.  For non-$p$--groups, other irreducible modules exist.  It is thus natural to generalize the Freyd conjecture to our presheaf context by letting $\sB$ be the category whose objects are the suspensions of the irreducible $kG$--modules.  The thick subcategory $\sC$ generated by this $\sB$ is then the subcategory of compact objects in the stable module category and one can ask when Conjecture \ref{FC} holds.   This presheaf version of the generating hypothesis in $\mathrm{StMod}_{kG}$ 
takes into account all the generating objects in the stable module category.
\end{exmp} 

\section{A general line of argument}

The presheaf perspective suggests a method of attack on the generalized Freyd conjecture. 
We have the forgetful functor 
\[ \bU = \io^*\colon \sP\sC\rtarr \sP\sB\]
given by restricting presheaves defined on $\sC^{op}$ to the full subcategory
$\sB^{op}$.  The functor $\bU$ has a left adjoint prolongation functor
\[ \bP = \io_!\colon \sP\sB\rtarr \sP\sC. \]
For $T\in \sP\sB$ and $K\in \sC$, $\bP T(K)$ is the categorical tensor product
\[ \bP T(K) = T\otimes_{\sB}\sC(K,-). \]
See, for example, \cite[I\S3]{MMSS} or \cite[I\S2]{MM}.  Since $\sB$ is a full
subcategory of $\sC$, the unit $\Id\rtarr \bU\bP$ of the adjunction is a natural
isomorphism \cite[I.3.2]{MMSS}.  We focus attention on the counit 
$\epz\colon \bP \bU\rtarr \Id$. 

Observe that the Freyd functor $\bF\colon \sT\rtarr \sP\sB$ is
the composite $\bU\bY$.  This leads to the following observation.  

\begin{prop}\label{crit1} 
The Freyd conjecture holds for  $(\sC,\sB)$ if 
$$ \epz\colon \bP\bF X = \bP\bU\bY X\rtarr \bY X $$ 
is an isomorphism
for all $X\in \sC$.
\end{prop}
\begin{proof}  Since the unit of the adjunction $(\bP,\bU)$ 
is an isomorphism, $\bU\epz$ is an isomorphism by one of the 
triangle identities.  Thus $\bF \iso \bU \bP\bF$.  Therefore, for a map 
$f\colon X\rtarr Y$ in $\sC$, $\bF f= 0$ if and only if $\bP\bF f =0$.  
By the Yoneda lemma, $f=0$ if and only if $\bY f =0$.   If $\epz$ is an 
isomorphism, then $\bP\bF f =0$ if and only if $\bY f =0$.
\end{proof}

In fact, less is needed.  Consider $\bY X(X) =[X,X]$. We will shortly
prove the following result. 

\begin{prop}\label{crit2} The Freyd conjecture holds for  $(\sC,\sB)$ if 
the identity map of $X$ is in the image of $\epz$ for all $X\in\sC$.
\end{prop}

We have the following starting point towards verification of the hypothesis.

\begin{lem}\label{start} The map
$\epz\colon \bP\bF X(J)\rtarr \bY X(J)$ 
is an isomorphism for all $X\in\sT$ and $J\in\sB$.
\end{lem}
\begin{proof} We have seen that $\bU\epz$ is an isomorphism,
and by definition $\bU T(J) = T(J)$ for any $J\in\sB$ and any
$T\in\sP\sC$.
\end{proof}

Now consider an exact triangle
\begin{equation}\label{cof}
\xymatrix@1{
K\ar[r] & L\ar[r] & M\ar[r] &  \SI K\\}
\end{equation} 
in $\sT$, where $K$, $L$, and $M$ are in $\sC$.   We have the commutative diagram

{\small

\begin{equation}\label{lad}
 \xymatrix{
\cdots \ar[r] & \bP\bF X (\SI K) \ar[r] \ar[d]_{\epz} 
& \bP\bF X (M) \ar[r] \ar[d]^{\epz} & \bP\bF X (L) \ar[r] \ar[d]^{\epz}
& \bP\bF X (K) \ar[r] \ar[d]^{\epz} & \cdots \\
\cdots \ar[r] & \bY X (\SI K) \ar[r] 
& \bY X (M) \ar[r]  & \bY X (L) \ar[r] 
& \bY X (K) \ar[r]  & \cdots. \\} 
\end{equation}

}

Since $\bY X (K) = [K,X]$, the lower row is exact.  By definition, $\sC$
is the smallest subcategory of $\sT$ that contains $\sB$, is closed under
retracts, and has the property that if two terms of an exact triangle are 
in $\sC$ then so is the third.  By an easy retract 
argument and the five lemma, this gives the following conclusion.

\begin{prop}\label{crit3} If the top row of Diagram (\ref{lad}) is exact 
for every exact triangle (\ref{cof}) and every $X\in\sC$, then 
$\epz\colon \bP\bF X \rtarr \bY X$ is an isomorphism for
every $X\in\sC$ and the Freyd conjecture holds for $(\sC,\sB)$.
\end{prop}

There is a reinterpretation of the exactness hypothesis that makes
it reminiscent of the standard result that the adjoint (if it exists) of an exact 
functor between triangulated categories is exact.  For $K$ and $X$
in $\sC$, the abelian group $\bP\bF X(K)$ is the coequalizer in $\sA\!b$
displayed in the diagram
\begin{equation}\label{PF} \xymatrix{
\sum_{I,J\in \sB} \sC(J,X)\otimes \sB(I,J)\otimes \sC(K,I)
\ar@<1ex>[d] \ar@<-1ex>[d] \\
\sum_{J\in\sB} \sC(J,X)\otimes \sC(K,J) \ar[d]\\
\sC(-,X)\otimes_{\sB}\sC(K,-),\\}
\end{equation}
where the parallel arrows are given by composition in $\sC$.
We use this to interpolate the proof of Proposition \ref{crit2}.  The composition maps
\[ \com\colon \sC(X,Y)\otimes \sC(J,X)\rtarr \sC(J,Y)\] 
induce a pairing 
\[ \com\colon \sC(X,Y)\otimes \bP \bF X(K)\rtarr \bP\bF Y(K) \]
such that $f\com z = \bP \bF f(z)$ for $z\in \bP\bF X(K)$ and
the following diagram commutes:
\[ \xymatrix{
\sC(X,Y)\otimes \bP \bF X(K) \ar[r]^{\id\otimes\epz} \ar[d]_{\com}
& \sC(X,Y)\otimes \bY X(K) \ar[d]^{\com}\\
\bP\bF Y(K) \ar[r]_{\epz} & \bY Y(K).\\} \]
Take $K=X$ and suppose that $\epz(z) = \id_X$ (for the top map $\epz$).   If $\bF f = 0$, then 
$\epz(f\com z) = \epz \bP\bF f(z) = 0$.  By the diagram, this equals
$f\com \epz(z) = f$ and so $f=0$.

Let us write $\sP_d\sB$ and $\sP_d\sC$ for the categories of covariant
additive functors on $\sB$ and $\sC$, and similarly write $\bU_d$,
$\bP_d$, $\bF_d$, and $\bY_d$ for the corresponding functors.  (The $d$ 
stands for dual.) We 
are just interchanging $\sB$ and $\sC$ with their opposite categories.
Visibly, we again have $\bU_d\bY_d = \bF_d$ and again have an adjunction
$(\bP_d,\bU_d)$ with $\bU_d\bP_d\iso\Id$.  By symmetry, we have
\begin{equation}\label{P'F'}
\bP\bF X(K) = \bP_d\bF_d K(X). 
\end{equation}
But in this dual reformulation, the exactness hypothesis on $K$ for fixed
$X$ is now a levelwise exactness statement about the composite functor 
$\bP_d\bF_d\colon \sC\rtarr \sP_d\sC$.  Since $\bF_dK(J) = \sC(K,J)$ for
$J\in\sB$, $\bF_d$ clearly takes exact triangles in the variable $K$
to exact sequences for each fixed $J$. Thus a more general question to ask is 
whether or not $\bP_d\colon \sP_d\sB\rtarr \sP_d\sC$ preserves levelwise exactness.  
That is, is it true that if $T'\rtarr T\rtarr T''$ is a sequence of diagrams 
$\sB\rtarr \sA\!b$ such that the sequence $T'(J)\rtarr T(J)\rtarr T''(J)$ is exact 
for each $J\in \sB$, then the sequence
$\bP T'(X)\rtarr \bP T(X)\rtarr \bP T''(X)$ is exact for all $X\in \sC$?

Observe that we have not yet used any hypothesis on $\sB$, other than
that it generates the thick subcategory $\sC$ of the triangulated category
$\sT$.  Thus all that we have done
is to give a purely formal reduction of the general problem. 

\section{The derived category of ring}

Our framework for understanding the Freyd conjecture leads to a transparent proof of 
the result of Lockridge \cite[3.9]{L} that the generalized Freyd conjecture holds in the derived category of a von Neumann regular ring.  
We simply observe that the hypotheses of Proposition \ref{crit3} hold in this case by one
of the equivalent definitions of a von Neumann regular ring.  However, our methods
do not prove the converse, which is proven in \cite{HLP, L}. 

Using
right $R$-modules for definiteness, let $\sD(R)$ be the derived category of a ring $R$ and let $\sB$ be the full subcategory of $\sD(R)$ whose objects are the shifts $\SI^i R$ of the chain complex that is $R$ concentrated in degree $0$. Then $\sC$ is the category of perfect chain complexes, namely those isomorphic in $\sD(R)$ to bounded chain complexes of finitely generated projective $R$-modules.
The Freyd functor assigns the homology groups $\sB(\SI^iR,X) = H_i X$ to a chain complex $X$, and Conjecture \ref{FC} is the assertion that a map between perfect chain complexes is $0$ in $\sD(R)$ if it induces
the zero map on homology.   Defining $H^i(K) = \sC(K,\Sigma^iR)$, as usual, we have the following observation in this case.

\begin{lem}\label{chaincxlem} For perfect chain complexes $K$ and $X$, $\bP\bF X(K)$ is isomorphic
to the abelian group $\sum_i H_i(X)\otimes_R H^i(K)$.
\end{lem}
\begin{proof}
By definition, $\bP\bF X(K)$ is the coequalizer displayed in the diagram
\[
\xymatrix{
\sum_{i,j\in \bZ} \sC(\Sigma^jR,X)\otimes \sB(\Sigma^iR,\Sigma^jR)\otimes \sC(K,\Sigma^iR)
\ar@<1ex>[d] \ar@<-1ex>[d] \\
\sum_{i\in\bZ} \sC(\Sigma^iR,X)\otimes \sC(K,\Sigma^iR) \ar[d]\\
\sC(-,X)\otimes_{\sB}\sC(K,-).\\}
\]
There are no maps $\Sigma^iR\to\Sigma^jR$ unless $i=j$, when $\sB(\Sigma^iR,\Sigma^iR)\cong R$.
The composition maps 
\[ \sC(\Sigma^iR,X)\otimes \sB(\Sigma^iR,\Sigma^iR)\rtarr \sC(\Sigma^iR,X) \]
specify the right action of $R$ on $H_i(X)$. Therefore our coequalizer diagram
can be rewritten as 
\[
\xymatrix{
\sum_{i\in \bZ} H_i(X)\otimes R\otimes H^i(K)
\ar@<1ex>[d] \ar@<-1ex>[d] \\
\sum_{i\in\bZ} H_i(X)\otimes H^i(K) \ar[d]\\
\sum_{i\in\bZ} H_i(X)\otimes_R H^i(K).}
\]
The conclusion follows.
\end{proof}

\begin{prop}\label{flathomology} If the homology $R$-modules $H_i(X)$ of any perfect chain 
complex $X$ are $R$-flat, then the Freyd conjecture holds for the pair $(\sC,\sB)$ in $\sD(R)$.
\end{prop}
\begin{proof}
Let $K\to L\to M\to \Sigma K$ be an exact triangle in $\sD(R)$. Since the functor 
$\sC(-,\Sigma^iR)$ takes exact triangles to exact sequences, the sequence 
\[H^{i-1}(K) \to H^i(M) \to H^i(L)\to H^i(K)\]
is exact.  By our flatness hypothesis, this sequence remains exact on tensoring with each
$H_i(X)$.  By Lemma \ref{chaincxlem}, when we take the direct sum over $i$ of these sequences, we obtain 
the exact sequence
\[\dotsb\rtarr \bP\bF X(\Sigma K)\rtarr \bP\bF X(M)\rtarr \bP\bF X(L)
\rtarr \bP\bF X(K)\rtarr \dotsb\]
The conclusion follows from Proposition 2.6. 
\end{proof}

By one definition, the ring $R$ is von Neumann regular if every $R$-module is flat.

\begin{cor} The Freyd conjecture holds for the derived category $\sD(R)$ of a 
von Neumann regular ring $R$.
\end{cor}

We have the following more general analogous result. 

\begin{cor} The Freyd conjecture holds for the derived category $\sD(R)$ if $R$ has weak 
dimension at most one and all finitely generated submodules of finitely generated projective
$R$-modules are FP-injective.
\end{cor}
The ring $R$ has weak dimension at most one if all submodules of flat modules are flat.  
A module $M$ is FP-injective if $\Ext^1(F,M)=0$ for all finitely presented modules $F$.  
By \cite[Theorem 4.89]{Lam}, $M$ is FP-injective if and only if any short exact sequence of 
right modules
\[0\to M\to L\to N\to 0\]
is pure exact, meaning that it remains exact on tensoring with any left module.
\begin{proof}
Let $R$ satisfy the stated hypotheses.  We show that the homology modules $H_i(X)$ 
of a perfect complex $X$ are flat.  Proposition \ref{flathomology} then applies.

Consider the short exact sequence
\begin{align}\label{SES}
0\to B_i(X)\to Z_i(X)\to H_i(X)\to 0.
\end{align}  The module $Z_i(X)$ is flat since it is a submodule of the projective, hence flat, module $X_i$.  By \cite[Corollary 4.86]{Lam}, it follows that $H_i(X)$ is flat if and only if 
the sequence (\ref{SES}) is pure exact.  The module $B_i(X)$ is finitely generated since it is a quotient of the finitely generated module $X_{i+1}$, and it is a submodule of the finitely 
generated projective $R$-module $X_i$.  By assumption, this implies that $B_i(X)$ is FP-injective.  Therefore (\ref{SES}) is pure exact and $H_i(X)$ is flat.
\end{proof}

This elementary result comes close to one direction of the best possible result about 
$\sD(R)$, which is proven in \cite[Theorem 2.1]{HLP} and
states that the Freyd conjecture holds for $\sD(R)$ if and only if $R$ has weak dimension at
most one and all finitely presented $R$-modules are FP-injective.


\end{document}